\documentclass[reqno, a4paper, 11pt]{amsart}
\usepackage{amssymb, amsthm, amsmath}
\numberwithin{equation}{section}

\renewcommand{\d}{\mathrm{d}}
\newcommand{\x}{\ma{x}}
\newcommand{\w}{\ma{w}}
\newcommand{\y}{\ma{y}}
\newcommand{\z}{\ma{z}}
\renewcommand{\v}{\ma{v}}
\renewcommand{\u}{\ma{u}}

\newcommand{\mbeta}{\boldsymbol{\beta}}

\newcommand{\R}{\mathbb{R}}
\newcommand{\F}{\mathbb{F}}

\newcommand{\Z}{\mathbb{Z}}
\newcommand{\N}{\mathbb{N}}
\newcommand{\Q}{\mathbb{Q}}

\newcommand{\bfP}{\mathbb{P}}

\newcommand{\ov}{\overline}
\newcommand{\ma}{\mathbf}

\renewcommand{\le}{\leqslant}
\renewcommand{\ge}{\geqslant}
\renewcommand{\leq}{\leqslant}
\renewcommand{\geq}{\geqslant}
\newcommand{\ben}{\begin{enumerate}}
\newcommand{\een}{\end{enumerate}}
\newcommand{\eit}{\begin{itemize}}
\newcommand{\beq}{\begin{equation}}
\newcommand{\eeq}{\end{equation}}
\newcommand{\ve}{\varepsilon}

\newcommand{\al}{\alpha}

\renewcommand{\mod}{\hspace{-0.2cm}\pmod}
\newcommand{\tmod}{\hspace{-0.1cm}\pmod}
\newcommand{\lab}{\label}

\newtheorem{thm}{Theorem}
\newtheorem*{thm*}{Theorem}
\newtheorem{lem}{Lemma}

\newtheorem*{cor*}{Corollary}

\DeclareMathOperator{\sing}{sing}

\theoremstyle{definition}

\newtheorem*{notat}{Notation}
\newcommand{\hcf}{\mathrm{gcd}}
\newcommand{\colt}[2]{\genfrac{}{}{0pt}{1}{#1}{#2}}

\begin{document}

\title{Integral Points on Cubic Hypersurfaces}

\author{T.D. Browning}

\author{D.R. Heath-Brown}

\address{School of Mathematics,  University of Bristol, Bristol BS8 1TW}
\email{t.d.browning@bristol.ac.uk}

\address{Mathematical Institute,
24--29 St. Giles', Oxford OX1 3LB}

\email{rhb@maths.ox.ac.uk}

\date{\today}

\maketitle


\section{Introduction}

Let $g\in\Z[x_1,\ldots,x_n]$ be an absolutely irreducible cubic
polynomial whose homogeneous part is non-degenerate. 
The primary goal of  this paper is to investigate the set
of integer solutions to the equation $g=0$.  Specifically, we 
shall try to determine
conditions on $g$ under which we can show that there 
are infinitely many solutions. An obvious necessary condition for the 
existence of integer solutions
is that the congruence
\begin{equation}
  \label{eq:cong_cond}
  g(x_1,\ldots,x_n)\equiv 0 \mod{p^k}, 
\end{equation}
should be soluble for every prime power $p^k$. We shall henceforth
refer to this condition as  ``the Congruence Condition''.

There is no condition on the
size of $n$ sufficient to ensure that \eqref{eq:cong_cond} is always 
soluble for
non-homogeneous $g$. In fact, even when the Congruence Condition is
satisfied for a polynomial $g$, and $n$ is large, this is still not 
sufficient to ensure
the existence of integer solutions to the corresponding equation $g=0$.  An
illustration of this is provided by the polynomial
\begin{equation}
  \label{eq:watson}
g=(2x_1-1)(1+x_1^2+\cdots +x_n^2)+x_1x_2,  
\end{equation}
discovered by Watson. 
Now it can be shown quite easily that $g$ satisfies the
Congruence Condition for $n\geq 2$. However, the equation $g=0$ is
insoluble in integers,  since $|2x_1-1|\geq 1
$ for every $x_1\in \Z$, so
$
|g|\geq 1+x_1^2+x_2^2-|x_1x_2|>0.
$

In view of Watson's example, it will be necessary to introduce an
auxiliary condition on the polynomials $g$ that we are able to handle.
Throughout this work we shall write $g_0$ for the homogeneous
cubic part of a polynomial $g\in\Z[x_1,\ldots,x_n]$. The
condition that we shall work with is phrased in terms of the
singular locus  of the projective hypersurface $g_0=0$ in
$\bfP_{\Q}^{n-1}$, which we denote by $\sing(g_0)$, 
a proper projective subvariety of $\bfP_{\Q}^{n-1}$. We set
\[  s(g_0):=\dim \sing(g_0)\]
for its projective dimension.  
Following the convention that $s(g_0)=-1$ if and
only if $g_0$ is a non-singular cubic form, we see that 
$s(g_0)$ is an integer contained in the interval
$[-1,n-2]$.
We are now ready to state our main result.

\begin{thm}\lab{main1}
Suppose that $g\in\Z[x_1,\ldots,x_n]$ is a cubic polynomial that
satisfies the Congruence Condition, such that $g_0$ is non-degenerate,
and having 
$$
n\geq 11+ s(g_0).
$$
Then the equation $g=0$ has infinitely many solutions in integers.
\end{thm}

This improves on work of Skinner \cite{skinner}, who has established
the same conclusion under the assumption that 
$n\geq 18+ s(g_0)$. At this point we note that the polynomial 
\eqref{eq:watson} has
homogeneous cubic part $g_0=2x_1(x_1^2+\cdots+x_n^2)$, which defines a
reducible cubic hypersurface with singular locus of dimension 
$
s(g_0)=n-3.
$
Hence Theorem \ref{main1} is not applicable to this particular example.

It is interesting to place Theorem \ref{main1} in the context of other
work in the literature. This topic has been extensively studied only in the
framework of homogeneous $g$. For arbitrary cubic forms the best result
available was, until very recently, due  
to Davenport \cite{D16}. This shows that there exists a non-trivial integer
solution to the homogeneous equation $g=0$ as soon as $n\geq 16$. This 
has now been
improved by the second author \cite{14}, who has shown that $n\geq 14$
variables are enough. One can do even better when the form under 
consideration is non-singular,
and the second author has shown that $n\geq 10$ variables suffice for
such forms \cite{hb-10}. This in turn has been sharpened by Hooley in 
a series of
papers \cite{hooley1, hooley2, hooley3}, to the extent that when 
$n\geq 9$ and the Congruence
Condition is satisfied, then 
the homogeneous equation $g=0$ is soluble in integers provided that the
corresponding hypersurface has only finitely many singularities and
these are linearly independent double points.

Returning to the topic of arbitrary cubic polynomials
$g\in\Z[x_1,\ldots,x_n]$, Davenport and Lewis \cite{DL} have also
considered the problem of determining when the equation $g=0$ has an
integer solution.  Their main results are phrased in terms of the
so-called $h$-invariant. Let $g_0$ denote the cubic homogeneous part
of $g$, as above. Then the invariant $h=h(g_0)$ is defined to be the
least positive integer such that $g_0$ can be written identically as
$$
g_0(x_1,\ldots,x_n)=L_1Q_1 +\cdots+L_hQ_h,
$$
for linear forms $L_i\in\Z[x_1,\ldots,x_n]$ and quadratic forms
$Q_i \in\Z[x_1,\ldots,x_n]$.
With this notation in mind, Davenport and Lewis show that the 
equation $g=0$ has infinitely many solutions in integers provided that 
$g$ satisfies the Congruence Condition, and has $h(g_0)\geq 17$.
In the course of generalising this work to the setting of arbitrary number
fields, Pleasants \cite{pleasants} has improved this lower bound to
$h(g_0)\geq 16$.
In a series of papers, culminating in \cite{watson3}, Watson
has shown 
that the equation $g=0$ is soluble in integers provided that 
$g$ satisfies the Congruence Condition, and has
$$
4\leq h(g_0)\leq n-3.
$$
One may combine this result with \cite{pleasants}, to conclude that 
the equation $g=0$ is soluble provided that 
$g$ satisfies the Congruence Condition, and has
$$
n\geq 18, \quad h(g_0)\geq 4.
$$
Note that one has $h(g_0)=1$ in \eqref{eq:watson}, so that none of
these results apply to this particular example.

Outside of the work of Skinner \cite{skinner}, it is not entirely
straightforward to compare the relative merits of Theorem \ref{main1}
with this previous body 
of work. It is true, however, that the condition on $s(g_0)$
is much easier to check than the condition on $h(g_0)$.
Theorem \ref{main1} also has something new to say about 
the case in which $g$ is homogeneous. In this setting one can view it as a bridge between \cite{hb-10} and 
\cite{14}, giving a new result for cubic forms $g$ such that
$$
13\geq n\geq 
\left\{
\begin{array}{ll}
11, & \mbox{if $s(g)=0$,}\\
12, & \mbox{if $s(g)=1$,}\\
13, & \mbox{if $s(g)=2$.}
\end{array}
\right.
$$
It is a natural question whether the approach of Hooley
\cite{hooley1, hooley2, hooley3} can be adapted to handle polynomials
rather than forms.  However we shall content ourselves with
investigating the extension of the second author's methods
\cite{hb-10}, since Hooley's technique is considerably more
delicate.  Our strategy will be to prove Theorem \ref{main1} for the case in
which $g_0$ is non-singular, that is to say that $s(g_0)=-1$,
and subsequently to deduce the general case via a hyperplane 
slicing argument.

Much of the work in this paper consists of trivial generalizations of
the second author's paper \cite{hb-10}.  However there are three new
things to be done.  Firstly, we have new complete exponential sums to
consider, which we shall reduce to Deligne's 
results \cite{deligne}, through a technique of Hooley \cite{hooleysum}.
Secondly, we must reconsider the singular series and the Congruence
Condition.  Thirdly, we shall show how to treat
polynomials for which $g_0$ is singular.

\begin{notat}
Throughout our work $\N$ will denote the set of positive
integers.  For any $\al\in \R$, we shall follow common convention and
write $e(\al):=e^{2\pi i\al}$ and $e_q(\al):=e^{2\pi i\al/q}$. 
All order constants will be allowed to depend on the polynomial $g$.
\end{notat}

\section{The circle method}\label{s:act}

In this section we recall the framework of the Hardy--Littlewood
circle method, as it applies to our problem on cubic polynomials. 
Our approach is based on the version of the circle method due to the
second author \cite{hb-10}, which incorporates Kloosterman's method for
tracking the precise location of the endpoints in the
decomposition of the unit interval into Farey arcs.  We have decided to
follow \cite{hb-10} as closely as possible, rather than to
incorporate some of the improvements introduced by Hooley.  We hope
readers will appreciate having to familiarize themselves with only one
source, rather than two.

We begin by choosing a non-zero real point $\x_0$ on $g_0(\x)=0$,
satisfying the additional condition that the matrix of second
derivatives of $g_0$ should have rank at least $n-1$ at $\x_0$.  The
existence of such a point is established as \cite[Lemma 5]{hb-10}. We
let $P$ be a large parameter, which we think of as tending to
infinity, and we define the weight
\[ w(\x):=\exp(-|\x-P\x_0|^2P_0^{-2}),\]
where $P_0:=P(\log P)^{-2}$.  Our main task is then to examine
the asymptotic behaviour of 
\begin{equation}
  \label{eq:Nw}
N(g;P):= \sum_{\colt{\x\in \Z^n}{g(\x)=0}} w(\x),  
\end{equation}
as $P\rightarrow \infty$. 
We define the singular series
\[
\mathfrak{S}(g):= \sum_{q=1}^\infty  
q^{-n}\sum_{\colt{a=1}{\hcf(a,q)=1}}^q \sum_{\x\mod{q}}e_q(ag(\x)), 
\]
when it converges, and the singular integral
\[\mathfrak{I}(g;P):=
\int_{-1}^{1} 
\int_{\R^n} w(\x)e(z g(\x))\d\x\d z.
\]
Then we shall prove the following estimate.

\begin{thm}\label{main1'}
Let $g\in  \Z[x_1,\ldots,x_n]$ be a cubic polynomial for which $g_0$
is non-singular. Assume that $n\geq 10$. Then there exists $\delta >0$ 
such that  
\begin{equation}
  \label{eq:125}
N(g;P)=\mathfrak{S}(g)\mathfrak{I}(g;P) + O\big(P^{n-3-\delta}\big).
\end{equation}
We have
\[
P^{n-3}(\log P)^{2-2n}\ll\mathfrak{I}(g;P)\ll P^{n-3}(\log P)^{2-2n}.
\]
Moreover $\mathfrak{S}(g)>0$ providing that the Congruence
Condition holds.
\end{thm}

As is implicit in the statement of Theorem \ref{main1'}, 
the singular series $\mathfrak{S}(g)$ is convergent for $n\geq 10$. In
fact we shall establish absolute convergence under this hypothesis.

Define the cubic exponential sum
\[
S(\al):=\sum_{\x\in\Z^n}w(\x)e(\al g(\x)),
\]
for $P\geq 2$. Then $S(\al)$ converges absolutely,
and for any $Q\geq 1$ we have 
$$
N(g;P)=\int_0^1 S(\al)\d\al=
\int_{-\frac{1}{1+Q}}^{1-\frac{1}{1+Q}}S(\al)\d\al,
$$
where $N(g;P)$ is given by \eqref{eq:Nw}.
In the form of the circle method developed by the second author
\cite{hb-10}, one proceeds to break the interval 
$[-\frac{1}{1+Q}, 1-\frac{1}{1+Q}]$ 
according to the Farey dissection of order $Q$. This ultimately yields
\begin{equation}
  \label{eq:Nw'}
\begin{split}
N(g;P)=&\sum_{q\leq Q}\int_{-\frac{1}{qQ}}^{\frac{1}{qQ}}S_0(q;z)\d
 z +O\big(Q^{-2}E(g;P,Q)\big),
\end{split}
\end{equation}
for any $Q\geq 1$, where 
\[
  \label{eq:error}
  E(g;P,Q):=
\sum_{q\leq Q}\sum_{|u|\leq
\frac{q}{2}}\frac{\max_{\frac{1}{2}\leq qQ|z|\leq 1}|S_u(q;z)|}{1+|u|},
\]
and 
\[
  \label{eq:Su}
S_u(q;z):=\sum_{\colt{a=1}{\hcf(a,q)=1}}^q e_q(\ov{a}u)S(a/q+z).
\]
This is \cite[Lemma 7]{hb-10}. 
We shall find that our work is optimised by taking $Q=P^{3/2}$ in
\eqref{eq:Nw'}.  As in \cite[\S 4]{hb-10} we shall estimate $S_u(q;z)$
via an application of the Poisson
summation formula.  This leads to the expression
$$
S_u(q;z)=q^{-n}\sum_{\v\in\Z^n}S_u(q;\v)I(z;q^{-1}\v),
$$
where
\[
  \label{eq:Sq}
S_u(q;\v):=\sum_{\colt{a=1}{\hcf(a,q)=1}}^q e_q(\ov{a}u)
\sum_{\y\mod{q}}e_q(ag(\y)-\v.\y),
\]
and 
\[
  \label{eq:Iq}
  I(z;\mbeta):=\int w(\x)e(z g(\x)+\mbeta.\x)\d\x.
\]
We proceed to estimate $I(z;\mbeta)$ as in \cite[\S 4]{hb-10}.  Some
small modifications are needed.

We begin by correcting a minor error
in the original treatment.  It is necessary that the parameter $l$
introduced just before \cite[(4.4)]{hb-10} should depend only on
$\x_1$, and not on $t_1,\ldots,t_{i-1},t_{i+1},\ldots,t_n$.  Thus one
should take 
\[l=z\frac{\partial F}{\partial x_i}(\x_1)+\beta,\]
rather than $l=zf'(x_1)+\beta$.  The appropriate estimate for
$f'(x_1+u)$ is then
\[f'(x_1+u)=\frac{\partial F}{\partial
  x_i}(\x_1)+O(P_1L|x_1|)+O(P_1^2L^2),\] 
where $L:=\log (P(2+|\mbeta|))$.

We have now to consider the measure of the set of vectors $\x_1$ for
which (4.7) and (4.8) of \cite{hb-10} hold, that is to say, for which
\[
|\x_1-P\x_0|\ll L P_0 
\]
and
\[
|z\nabla g(\x_1)+\mbeta|\ll L^7(P^{-1}+|z|^{1/2}P^{1/2}). 
\]
Since $\nabla g(\x_1)=\nabla g_0(\x_1)+O(P)$, this latter constraint
yields
\[
|z\nabla g_0(\x_1)+\mbeta|\ll L^7 
(P^{-1}+|z|^{1/2}P^{1/2}+|z|P). 
\]
We can then proceed exactly as before so as to
deduce the following extension of \cite[Lemma 8]{hb-10}.

\begin{lem}\label{l8}
For $|z|\le 1$ we have
\[S_u(q;z)=q^{-n}\sum_{\colt{\v\in\Z^n}{|\v|\leq V_0}} 
S_u(q;\v)I(z;q^{-1}\v)+O(1),\]
with
\[V_0\ll (\log P)^7 q(P^{-1}+|z|P^2).\] 
Moreover for $|z|\le 1$ we have
\[ I(z;\mbeta)\ll (\log
P)^{7n}\big(P+\min\{P^n\,,\,P^{(3-n)/2}|z|^{(1-n)/2}\}\big). 
\]
\end{lem}

\section{The sum $S_u(q;\v)$ when $q$ is prime}\label{sec:katz} 

The sum $S_u(q;\v)$ satisfies the multiplicativity property 
\[
S_u(rs;\v)=S_{\bar{s}^2u}(r;\bar{s}\v)S_{\bar{r}^2u}(s;\bar{r}\v),\quad \mbox{for $\hcf(r,s)=1$},
\]
where $\bar{r},\bar{s}$ are any integers such that $r\bar{r}+s\bar{s}=1$.
Note that this reduces to (5.1) in \cite{hb-10} when $g$ is homogeneous.
It therefore suffices to examine the case in which $q$ is a prime
power.  In this section we handle prime values of $q$, 
using the following lemma, which summarizes the technique developed by Hooley
\cite{hooleysum}.

\begin{lem}\label{genhoo}
Let $F$ and $G$ be polynomials over $\Z$, of degree at most $d$, and let
$$
S:=\sum_{\x\in\mathbb{F}_p^n, G(\x)=0}e_p(F(\x)), 
$$ 
for any prime $p$. For each $j\ge 1$ write
$$
N_j(\tau)=\#\{\x\in\mathbb{F}_{p^j}^n:\,G(\x)=0,\,F(\x)=\tau\},
$$
and suppose that for each $j$ there is a real 
number $N(j)$ such that
\begin{equation}\label{hoo1}
\sum_{\tau\in\mathbb{F}_{p^j}}|N_j(\tau)-N(j)|^2\ll_{d,n} p^{kj},
\end{equation}
where $k$ is an integer independent of $j$.  Then
$$
S\ll_{d,n}p^{k/2}.
$$
\end{lem}

Our second key tool is an extension of the famous result of 
Deligne \cite{deligne}, due to Hooley \cite{JNT}.

\begin{lem}\label{lemdel}
Let $q=p^j$ and let $H(x_1,\ldots,x_m)$ be a form of degree $d$ defined
over $\F_q$.  Assume that $p\nmid d$, and write $s\ge -1$ for the dimension
of the singular locus of $H=0$ in 
$\mathbb{P}^{m-1}(\ov{\F_q})$.  
Then 
\[\#\{\x\in\F_q^m: H(\x)=0\}=q^{m-1}+O_{d,m}(q^{(m+1+s)/2}).\]
\end{lem}

Note that if $H$ is not absolutely irreducible we interpret $s$ as the
dimension of the variety $\nabla H(\x)=\ma{0}$, so that $s\ge m-3$. 
The  result is therefore trivial for such $H$.  Indeed it remains true
even if $H$ vanishes identically, since then $s=m-1$.

Our basic estimate for $S_u(p;\v)$ is now provided by the following result.

\begin{lem}\label{katz}
Let $g\in \Z[x_1,\ldots,x_n]$ be a cubic
polynomial, and let $p$ be a prime. Suppose that
$g_0$ is non-singular modulo $p$, and that $p\nmid u$.
Then there is a constant $C(n)$ such that
$$
|S_u(p;\v)| \le C(n) p^{(n+1)/2}.
$$
\end{lem}

A more general result has been given recently by Katz \cite{katzweb},
but we shall give a shorter self-contained treatment, based on  
Lemma \ref{genhoo}. 

Lemma \ref{katz} is trivial for $p\le 3$, so we shall assume that $p\ge 5$. 
For the proof we set 
\[F(a,b,\x)=ub+ag(\x)-\v.\x,\;\;\; G(a,b,\x)=ab-1,\]
so that $S=S_u(p;\v)$ in the notation of Lemma \ref{genhoo}.  Then 
\[N_j(\tau)=\#\{(b,\x)\in\F_q^{n+1}:b^2u+g(\x)-b\v.\x-b\tau=0,\,b\not=0\},\]
where $q=p^j$.  We convert this into a question about projective
varieties by defining $\tilde{g}(z,\x):=z^3g(z^{-1}\x)$ and 
\[
H_{\tau}(b,z,\x):=ub^2z+\tilde{g}(z,\x)-bz\v.\x-bz^2\tau, 
\]
so that
\[
N_j(\tau)=\frac{1}{q-1}\#\{(b,z,\x)\in\F_q^{n+2}:
H_{\tau}(b,z,\x)=0,\,bz\not=0\}.
\]
Solutions with $z=0$ have $\tilde{g}(0,\x)=g_0(\x)=0$.  According to 
Lemma \ref{lemdel} one has $g_0(\x)=0$ for $q^{n-1}+O_n(q^{n/2})$ values
of $\x$, whence
\begin{align*}
\#\{(b,z,\x)\in\F_q^{n+2}: H_{\tau}(b,z,\x)=0,\,b &\neq 0,\,z=0\}\\
&=(q-1)q^{n-1}+O_n(q^{(n+2)/2}).
\end{align*}
When $b=0$ the equation $H_{\tau}=0$ reduces to $\tilde{g}(z,\x)=0$.  Let $V$
denote, temporarily, the singular locus of the variety defined by
$\tilde{g}(z,\x)=0$, and let $H$ denote the hyperplane given
by $z=0$.  Any points $(0,\x)$ on $V\cap H$ would satisfy
$\nabla g_0(\x)=\ma{0}$.
Since $g_0$ is assumed to be non-singular modulo $p$ we conclude that
$V\cap H$ is empty.  Thus $V$ has dimension at most  zero.  An
application of Lemma \ref{lemdel} now shows that
\[\#\{(b,z,\x)\in\F_q^{n+2}:H_{\tau}(b,z,\x)=0,\,b=0\}=q^n+O_n(q^{(n+2)/2}).\]
Finally, a third application of Lemma \ref{lemdel} yields
\[\#\{(b,z,\x)\in\F_q^{n+2}:H_{\tau}(b,z,\x)=0\}=q^{n+1}+O_n(q^{(n+3+s)/2}),\]
where $s$ is the dimension of the singular locus of the variety
$W_{\tau}\subseteq \mathbb{P}^{n+1}(\ov{\F_q})$
defined by $H_{\tau}(b,z,\x)=0$.

On combining our various results we conclude that
\[N_j(\tau)=(q-1)q^{n-1}+O_n(q^{(n+1+s)/2}).\]
We claim that $s=-1$ for all but $O_n(1)$ values of $\tau$ in
the algebraic closure $\ov{\F_q}$, 
and that $s=0$ for the remaining
values.  If we take $N(j)=(p^j-1)p^{nj-j}$ 
in Lemma \ref{genhoo} we
will then obtain the required estimate \eqref{hoo1} with $k=n+1$,
whence $S=S_u(p;\v)=O_n(p^{(n+1)/2})$ as required.

Taking partial derivatives, we see that singular points on $W_{\tau}$ satisfy
\begin{equation}\label{s1}
2ubz-z\v.\x-z^2\tau=0,
\end{equation}
\begin{equation}\label{s2}
ub^2+\frac{\partial \tilde{g}}{\partial z}(z,\x)-b\v.\x-2bz\tau=0
\end{equation}
and
\begin{equation}\label{s3}
\frac{\partial \tilde{g}}{\partial x_i}(z,\x)-bzv_i=0,\;\;\;(1\le i\le n).
\end{equation}
If $z=0$ then \eqref{s3} reduces to $\nabla g_0(\x)=\ma{0}$, whence
$\x=\ma{0}$, since $g_0$ is non-singular modulo $p$.  We then have
$b=0$ from \eqref{s2}.  Thus there can be no singular points with
$z=0$, so that \eqref{s1} may be replaced by
\begin{equation}\label{s11}
2ub-\v.\x-z\tau=0.
\end{equation}
We now eliminate $\tau$ from \eqref{s2} and \eqref{s11} to produce
\begin{equation}\label{s22}
-3ub^2+\frac{\partial \tilde{g}}{\partial z}(z,\x)+b\v.\x=0.
\end{equation}
It follows that all singular points on $W_{\tau}$, regardless of the value
of $\tau$, lie on the variety $U$, say, given by \eqref{s3} and
\eqref{s22}.  We proceed to examine the intersection $I$, say, of $U$  
with the hyperplane $z=0$.  For points on $I$ the equation \eqref{s3}
reduces to $\nabla g_0(\x)=\ma{0}$, whence $\x=\ma{0}$.  Since
$p\nmid u$ it then follows from \eqref{s22} that $b=0$, so that $I$ is
empty, as a subset of $\mathbb{P}^{n+1}(\ov{\F_q})$.
We therefore conclude that $U$ 
has at most dimension zero, and therefore contains $O_n(1)$ points $(b,z,\x)$.
Finally we conclude that the various varieties $W_{\tau}=0$ have
between them at most $O_n(1)$ singular points.  Since one has
$z\not=0$ for any singular point, as noted above, any singular point
determines exactly one corresponding value of $\tau$, via
\eqref{s11}.  It therefore follows that there are $O_n(1)$ values of
$\tau$ for which $W_{\tau}$ is singular, and that if $W_{\tau}$ is
singular then it has $s=0$.  This establishes the claim above, and
thereby completes the proof of Lemma \ref{katz}.

For the case in which $p\mid u$ the situation is more complicated.  We
consider the projective variety defined by $g_0(\x)=0$.  Then the 
dual variety is a hypersurface,
defined by an equation $g^*(\x)=0$, say.  We
now have the following estimate.

\begin{lem}\label{hooley}
Let $g\in \Z[x_1,\ldots,x_n]$ be a cubic
polynomial, and let $p$ be a prime. Suppose that
$g_0$ is non-singular modulo $p$.  Then there is a constant $C(n)$ such that
$$
|S_0(p;\v)| \le C(n) p^{(n+1)/2}(p,g^*(\v))^{1/2}.
$$
\end{lem}

As before this is trivial for $p\le 3$.  When $p\mid g^*(\v)$ we apply
the estimate
\[\sum_{\x\mod{p}}e_p(f(\x))\ll_{d,n} p^{n/2},\]
of Deligne \cite{deligne}, which applies to any polynomial $f$ over
$\F_p$ of degree $d$, in $n$ variables, whose homogeneous part is
non-singular modulo $p$.  Taking $f(\x)= ag(\x)-\v.\x$ we see that
\[
\sum_{\x\mod{p}}e_p(ag(\x)-\v.\x)\ll_n p^{n/2} 
\]
for $p\nmid a$.  Summing over $a$ yields a satisfactory bound when 
$p\mid g^*(\v)$.

For the general case we begin by
observing that $p\nmid \v$, since $p\nmid g^*(\v)$.
It follows that
$$
\sum_{\x \mod{p}}e_p(\v.\x)=0,
$$
whence
$$
S_0(p;\v)=\sum_{a,\x \mod{p}}e_p(ag(\x)-\v.\x)
=p\sum_{\colt{\x \mod{p}}{p\mid g(\x)}}e_p(-\v.\x)=pS,
$$
say.  It is possible to handle this by an application of Hooley's
method.  However a more general result due to Katz \cite{katzunp} is
already available.  To put our sum into the correct form for Katz' estimate,
we define $\tilde{g}(z,\y)=z^3g(z^{-1}\y)$ and substitute 
$\x\equiv z^{-1}\y\tmod{p}$.  If we then let $z$ run
over the residue classes coprime to $p$ we find that
\[S=\frac{1}{p-1}\sum_{z,\y}e_p(-z^{-1}\v.\y)\]
where $z,\y$ run over solutions of $\tilde{g}(z,\y)\equiv 0\tmod{p}$
with $z\not\equiv 0\tmod{p}$.  Thus we have
\[S=\sum_{(z,\y)\in V}e_p(-z^{-1}\v.\y),\]
where $V$ is the projective variety over $\F_p$ given by $\tilde{g}=0$
and $z\not=0$.  For this type of sum Katz \cite{katzunp}
shows that $S\ll_n p^{m/2}$, 
where $m=n-1$ is the dimension of $V$
in projective space, under the conditions that $\tilde{g}$ is
absolutely irreducible over $\F_p$, and that the variety
$\tilde{g}(z,\y)=\v.\y=z=0$ is smooth and of dimension $m-2=n-3$ in 
$\mathbb{P}^n(\ov{\F_p})$. Since $\tilde{g}(z,\y)=\v.\y=z=0$ implies
$g_0(\y)=\v.\y=0$, this second condition follows from our assumption
that $p\nmid g^*(\v)$.  Moreover $g_0$ is absolutely irreducible,
since it is non-singular, and the absolute irreducibility of
$\tilde{g}$ follows.  This completes our treatment of Lemma \ref{hooley}

\section{The sum $S_u(q;\v)$ when $q$ is square-full}\label{sec:sqf}

When $q$ is square-full we follow the analysis of \cite[\S 6]{hb-10},
with only minor modifications.  The sum $S_{\ma{k},\ma{h}}$ becomes
\[
S_{\ma{k},\ma{h}}=\sum_{\ma{j}\mod{q_2}} 
e_{q_2}(\ma{k}.\ma{j}+\frac{s}{2}\ma{j}^T\hspace{-1mm}M(\ma{h})\ma{j}),\]
where $M(\ma{h})$ is the matrix of second derivatives of $g(\ma{h})$.
Similarly $N(q_3;\ma{h})$ becomes
\[\tilde{N}(q_3;\ma{h})=
\#\{\ma{j}\mod{q_3}: q_3\mid\frac{1}{6}M(\ma{h})\ma{j}\}.\] 
We now have the following analogue of \cite[Lemma 4]{hb-10}.

\begin{lem}\label{anaL4}
Let
\[\tilde{N}(q):=
\#\{\ma{h},\ma{j}\mod{q}: q\mid\frac{1}{6}M(\ma{h})\ma{j}\}.\]
Then there is a constant $A$ such that
\[\tilde{N}(q)\le A^{\omega(q)}q^n\]
for every square-free $q$.  
\end{lem}

To prove this we observe that
\[M(\ma{h})\ma{j}=M_0(\ma{h})\ma{j}+M_1\ma{j},\]
where $M_0(\ma{h})$ is the matrix of second derivatives of $g_0$, and
$M_1$ is the matrix of second derivatives of the quadratic part of
$g$.  It follows that
\begin{align*}
\tilde{N}(q)&=\sum_{\ma{j}\mod{q}}
\#\{\ma{h}\mod{q}: q\mid\frac{1}{6}\big(M_0(\ma{h})\ma{j}+M_1\ma{j}\big)\}\\
&=\sum_{\ma{j}\mod{q}}
\#\{\ma{h}\mod{q}: q\mid\frac{1}{6}\big(M_0(\ma{j})\ma{h}+M_1\ma{j}\big)\}.
\end{align*}
However if $M_0$ is a square matrix and $\ma{c}$ is a constant vector
we claim that
\[\#\{\ma{h}\mod{q}: q\mid M_0\ma{h}+\ma{c}\}\le
\#\{\ma{h}\mod{q}: q\mid M_0\ma{h}\}.\]
It will then follow that
\begin{align*}
\tilde{N}(q)&\le\sum_{\ma{j}\mod{q}}
\#\{\ma{h}\mod{q}: q\mid\frac{1}{6}M_0(\ma{j})\ma{h}\}\\
&= \#\{\ma{j},\ma{h}\mod{q}: q\mid\frac{1}{6}M_0(\ma{j})\ma{h}\}.
\end{align*}
The lemma therefore follows from \cite[Lemma 4]{hb-10}, since the
final expression is just $N(q)$ for the non-singular form $g_0$.

It remains to prove the claim above.  If there is no vector $\ma{h}$
with $q\mid M_0\ma{h}+\ma{c}$ the result is trivial.  Otherwise let
$\ma{h}_0$ be any such vector.  Then $q\mid M_0\ma{h}+\ma{c}$ if and
only if $q\mid M_0(\ma{h}-\ma{h}_0)$, and the required bound follows.
This completes the proof of Lemma \ref{anaL4}.

We may now continue with the analysis as in \cite[\S 6]{hb-10},
finding, in the notation of \cite[page 242]{hb-10}, that
\[|S_1|^2\le\tilde{N}(q_3)|S_0|\ll \tilde{N}(q_3)q_1^nT(\ma{a}),\]
say, where
\[S_0=\sum_{\ma{h}_3\mod{q_4}}
e_{q_4}(q_3s\ma{a}.\nabla g_0(\ma{h}_3)+s\ma{a}^T\hspace{-1mm}M_1\ma{h}_3)
\sum_{\ma{h}_2\mod{q_1}}e_{q_4}(s\ma{a}^T\hspace{-1mm}M_0(\ma{h}_3)\ma{h}_2)\]
with
\[T(\ma{a})=\#\{\ma{h}_3\mod{q_5}:q_5\mid M_0(\ma{a})\ma{h}_3\}.\]

Everything now proceeds as before, leading to the following 
variant of \cite[Lemma 14]{hb-10}.  

\begin{lem}\label{anaL14}
There is a positive constant $A$,
such that for any integer vector $\v_0$ we have 
\[\sum_{|\v-\v_0|\le V}|S_u(q;\v)|\le   
A^{\omega(q)}(\log(q+1))^{2n}q^{n/2+1}(V^n+q^{n/3}),\] 
uniformly in $\v_0$, whenever $q$ is square-full. 
\end{lem}

The effect of introducing $\v_0$ into the analysis of \cite[\S  
6]{hb-10} is to modify the sum $S_2$ which occurs there. 
However the same estimate for $S_2$ still holds, and the proof 
goes through as before.

We shall also want to consider the sum in Lemma \ref{anaL14} with
$u=0$ and with $\v$
restricted by the condition $g^*(\v)=0$.  Here we follow the analysis
in \cite[\S 7]{hb-10}, with $\ma{h}.\nabla F(\ma{j})$ replaced by
$\frac{1}{2}\ma{j}^T\hspace{-1mm}M(\ma{h})\ma{j}$ throughout, so that
$N(q_3;\ma{h})$ becomes $\tilde{N}(q_3;\ma{h})$.  With these trivial
changes everything goes through as before, up to the 
treatment of the sum $S(q)$ defined in \cite[(7.3)]{hb-10}.
Let 
$$
q_1=\prod_{p^u\| q}p^{[u/2]},\quad 
q_2=\prod_{\colt{p^u\| q}{2\nmid u} }p, \quad  \mbox{and} \quad 
q_4=\prod_{\colt{p^u\| q}{2\nmid u, u\geq 13} }p.
$$
In the present setting we will only be able to establish that 
if $q$ is square-full then 
\begin{equation}
  \label{eq:SQ}
  S(q)\ll q_1^{n-1+\ve}q_2^{1/2}q_4^{1/2}.
\end{equation}
when $q$ is square-full. The corresponding bound in \cite{hb-10} 
is somewhat sharper, in that the factor $q_4^{1/2}$ is absent.    
We  shall prove (\ref{eq:SQ}) in a moment, but first we show how
it suffices for our purposes. 
Inserting (\ref{eq:SQ}) into \cite[(7.3)]{hb-10}, 
an application of \cite[Lemma 15]{hb-10} reveals that  
\[
\sum_{\colt{|\v|\le V}{g^*(\v)=0}}|S_0(q;\v)|\ll 
q_1^{2n+1+\ve}q_2^{(n+3)/2}q_4^{1/2} 
\big(1+Vq_1^{-1}\big)^{n-3/2}\log(V+1)
\]
whenever $q$ is square-full.
In order to establish the analogue of \cite[Lemma 16]{hb-10}
we are now left with a parallel calculation to the three lines at the
bottom of \cite[page 245]{hb-10}. Note that $q_4\leq q_2\leq q_1$ and
$q=q_1^2q_2$. Hence 
$$
q_1^{2n+1}q_2^{(n+3)/2}q_4^{1/2} \leq
q_1^{2n+1}q_2^{n/2+2} \leq q^{n+1/2}.
$$
Similarly, we have 
$$
q_1^{n+5/2}q_2^{(n+3)/2}q_4^{1/2} \leq
(q_1^2q_2)^{(2n+5)/4}q_2^{1/4}q_4^{1/2}.
$$
Thus it remains to confirm that $q_2^{1/4}q_4^{1/2}\leq q^{1/12}$, which
it suffices to verify at each prime power $q=p^e$. This is trivial when $e$ is even
since then $q_2=q_4=1$.  When $e=2f+1$ is odd we have
$$
q_2^{1/4}q_4^{1/2} =\left\{
\begin{array}{ll}
p^{1/4}, & \mbox{if $f \leq 5$, }\\
p^{3/4}, & \mbox{if $f \geq 6$, }
\end{array}
\right.
$$
which is always at most $p^{(2f+1)/12}$. 
Assuming the validity of \eqref{eq:SQ}, this therefore 
yields the following result, corresponding to \cite[Lemma 16]{hb-10}.

\begin{lem}\label{anaL16}
We have
\[
\sum_{\colt{|\v|\le V}{g^*(\v)=0}}|S_0(q;\v)|\ll
q^\ve\big(q^{n+1/2}+V^{n-3/2}q^{n/2+4/3}\big)\log(V+1),
\]
whenever $q$ is square-full.
\end{lem}

It remains to establish \eqref{eq:SQ}. Let $\v_0 \in \Z^n$ be an
arbitrary vector. Taking $V=q^{1/3}$, it follows from an application
of Lemma \ref{anaL14} that
\begin{equation}
  \label{eq:v0}
|S_0(q;\v_0)|\leq 
\sum_{|\v-\v_0|\le V}|S_0(q;\v)| \ll
q^{5n/6+1+\ve},  
\end{equation}
when $q$ is square-full.  Furthermore, the implied constant 
in this estimate does
not depend on $\v_0$. 
Arguing as in \cite{hb-10} one
easily checks that $S(q)$ is multiplicative in $q$, whence it suffices to
estimate $S(p^e)$ for $e\geq 2$. When $e$ is even, so that $q_3=1$,
the argument based on exponential sums in \cite[page 245]{hb-10} goes
through with no changes. This yields 
$$
S(p^{2f})=S_0(p^f) \ll p^{f(n-1)},
$$
for any $f\geq 1$. Indeed, once combined with Lemma \ref{katz}, the
bound in \eqref{eq:v0} gives $S_0(p^g;\ma{0})\ll
p^{g(5n/6+1+\ve)}$ for any $g\geq 1$. 
This argument also takes care of the finitely many primes $p$ 
for which $p\mid 6$ or for which 
the reduction of $g_0$ modulo $p$ is singular, since in 
these cases we have $\tilde{N}(p;\ma{h})\ll 1$. 

Turning to the case of odd $e\geq 2$, we suppose that $e=2f+1$ and
that $g_0$ is non-singular modulo $p$, with $p \nmid 6$. Thus $q_1=p^f$ and
$q_3=p$, and it follows that 
$$
S(p^e)=\sum_{\colt{\mathbf{k}\mod p}{p\mid g(\ma{k})}} 
\tilde{N}(p;\ma{k})^{1/2}M(p^f;\ma{k}),
$$
where
$$
M(p^f;\ma{k})= \#\big\{ \ma{h}\mod{p^{f}}: p\mid \ma{h}-\ma{k},  
~p^f\mid g(\ma{h})\big\}.
$$
When $p\nmid \nabla g(\ma{k})$ a straightforward argument 
based on Hensel's lemma
reveals that $M(p^f;\ma{k})\ll p^{(f-1)(n-1)}$. Thus the overall
contribution from such $\ma{k}$ is 
\begin{align*}
\ll p^{(f-1)(n-1)} \sum_{\colt{\mathbf{k}\mod p}{p\mid g(\ma{k})}} 
\tilde{N}(p;\ma{k})^{1/2} 
&\ll p^{f(n-1)+1/2},
\end{align*}
by an application of Cauchy's inequality and Lemma \ref{anaL4}.
To handle the contribution from the remaining $\ma{k}$, we observe
that there can only be $O(1)$ values of $\ma{k}$ modulo $p$ for which
$p\mid g(\ma{k})$ and $p\mid \nabla g(\ma{k})$. This follows
from the fact that the corresponding projectivised
variety has dimension $0$, as we saw in our analysis of $V$ in the proof of
Lemma \ref{katz}. Taking $\tilde{N}(p;\ma{k})^{1/2}\leq p^{n/2}$, and
incorporating our work above, we deduce that 
$$
S(p^e)\ll p^{f(n-1)+1/2}+ p^{n/2}
\max_{\mathbf{k}}
M(p^f;\ma{k})
$$
for $e=2f+1$, 
where the maximum is taken over all $\mathbf{k}$ modulo $p$ such that 
$p\mid g(\ma{k})$ and $p\mid \nabla g(\ma{k})$.
We will show that
\begin{equation}
  \label{eq:show}
  \max_{\mathbf{k}}
M(p^f;\ma{k})
\ll p^{f(n-1+\ve)+5-n+\theta_{p}(e)},
\end{equation}
where $\theta_p(e)=1$ if $e=2f+1$ with $f\geq 6$, and $\theta_p(e)=0$ otherwise. 
Once inserted into our bound for $S(p^e)$ this implies that
$$
S(p^e)\ll p^{f(n-1+\ve)}(p^{1/2}+ p^{5-n/2+\theta_p(e)}) \ll 
p^{f(n-1+\ve)+1/2+\theta_p(e)/2},
$$
since $n \geq 10$. In view of the fact that $\theta_p(e)=0$ unless
$e=2f+1\geq 13$, this is therefore enough to complete the proof of
\eqref{eq:SQ}.

We will use exponential sums to estimate $M(p^f)=M(p^f;\ma{k})$. Thus we find
that
\begin{align*}
M(p^f)
&=\frac{1}{p^{f+n}} \sum_{s\mod {p^f}}\sum_{\ma{j}\mod p}
\sum_{\ma{h}\mod {p^f}} e\big(s g(\ma{h})/p^f+\ma{j}.(\ma{h}-\ma{k})/p\big)\\
&=\frac{1}{p^{f+n}}\sum_{0\leq g\leq f} \sum_{\colt{t\mod {p^g}}{p\nmid t}}\sum_{\ma{j}\mod p}
\sum_{\ma{h}\mod {p^f}} e_{p^g}\big(t g(\ma{h})+p^{g-1}\ma{j}.(\ma{h}-\ma{k})\big),
\end{align*}
on splitting $s$ according to the value of the highest common factor $p^{f-g}$
of $s$ with $p^f$. Fix a choice of $\ell$, with $1\leq \ell \leq f$. 
Let us write $M_1(p^f)$ for the contribution  to $M(p^f)$ from values
of $g\leq \ell$, and $M_2(p^f)$ for the corresponding contribution from values
of $g>\ell$. Beginning with small values of $g$, we 
reverse the process above to deduce that 
\begin{align*}
M_1(p^f)
&=\frac{1}{p^{f}}\sum_{0\leq g\leq \ell} \sum_{\colt{t\mod {p^g}}{p\nmid t}}
\sum_{\colt{\ma{h}\mod {p^f}}{p\mid (\ma{h}-\ma{k})}} 
e_{p^g}\big(t g(\ma{h})\big)\\
&=\frac{1}{p^{f}}\sum_{0\leq g\leq \ell} p^{fn}  
\big(p^{(1-n)g} M(p^g) -  p^{(1-n)(g-1)} M(p^{g-1})\big)\\  
&=p^{f(n-1)} p^{\ell(1-n)} M(p^\ell).
\end{align*}
Here we have followed the convention that $M(p^{g-1})=0$ when $g=0$.
Employing the crude upper bound $M(p^{\ell})\ll p^{(\ell-1)n}$, we
deduce that
\begin{equation}\label{eq:M1}
M_1(p^f)
\ll p^{f(n-1)+\ell-n}.
\end{equation}

To produce a bound for $M_2(p^f)$, we apply \eqref{eq:v0} 
to deduce that 
\begin{align*}
\sum_{\colt{t\mod{p^g}}{p\nmid t}}\sum_{\ma{h}\mod {p^f}} e_{p^g}\big(t
g(\ma{h})+p^{g-1}\ma{j}.\ma{h}\big)
&=p^{(f-g)n} S_0(p^{g};-p^{g-1}\ma{j})\\
&\ll  p^{fn+g(1-n/6+\ve)},
\end{align*}
for each $g$ and $\ma{j}$.  Hence  
\begin{align*}
M_2(p^f)
&\ll\frac{1}{p^{f+n}}\sum_{\ell< g\leq f} \sum_{\ma{j}\mod p}
p^{fn+g(1-n/6+\ve)}\\
&\ll p^{f(n-1)} \sum_{\ell< g\leq f} 
 p^{g(1-n/6+\ve)}\\
&\ll p^{f(n-1+\ve)-\ell(n/6-1)}.
\end{align*}
When $f\geq 6$, we take $\ell=6$ and combine the above estimate with
\eqref{eq:M1} to conclude that
$$
M(p^f)=M_1(p^f)+M_2(p^f) \ll  p^{f(n-1+\ve)+6-n}= p^{f(n-1+\ve)+5-n+\theta_p(e)}.
$$
This is satisfactory for \eqref{eq:show}. Alternatively, when $f\leq
5$ we apply \eqref{eq:M1} with $\ell=f$ to deduce that
$$
M(p^f)=M_1(p^f) \ll  p^{f(n-1+\ve)+5-n}.
$$
This too is satisfactory for \eqref{eq:show}, and so completes its proof.

\section{Proof of Theorem \ref{main1'}}\label{endgame}

We now estimate the various terms in \eqref{eq:Nw'}, just as in
\cite[\S\S 8, 9 and 10]{hb-10}, finding that \eqref{eq:125} holds with
\begin{align*}
\mathfrak{S}(g)&=\sum_{q=1}^{\infty}q^{-n}S_0(q;\ma{0}),\\
\mathfrak{I}(g;P)& =\int_{-1}^{1}I(z;\ma{0})\d z, 
\end{align*}
and any fixed $\delta<1/18$.  
Note that the only difference arises from the additional term
$P (\log P)^{7n}$ appearing in the estimate for $I(z;\mbeta)$ in Lemma~\ref{l8}.
It is a trivial matter to check that this term makes a satisfactory
overall contribution to $N(g;P)$.
It therefore remains to show that
$\mathfrak{S}(g)$ is strictly positive, and that
\begin{equation}
  \label{eq:n-3}
P^{n-3}(\log P)^{2-2n}\ll\mathfrak{I}(g;P)\ll P^{n-3}(\log P)^{2-2n}.
\end{equation}

For the first task it suffices, as usual, to show that $g(\x)=0$ has a
non-singular $p$-adic integer zero for every prime $p$.  This problem
is discussed by Davenport and Lewis \cite[\S 2]{DL}, where it is shown
that if the Congruence Condition holds, then non-singular $p$-adic
solutions exist, except possibly in Case 3b of \cite[\S 2]{DL}.  In
the excluded case there is a non-singular matrix $M$ and a positive integer
$r\le 4$ such that
\[g(M\x)=x_1R_1(x_1,\ldots,x_n)+\cdots+x_rR_r(x_1,\ldots,x_n)+
R(x_1,\ldots,x_r)\]
for certain quadratic forms $R_1,\ldots,R_r$ and $R$.  If $g$ has such
a representation then 
\[g_0(M\x)=x_1R_1(x_1,\ldots,x_n)+\cdots+x_rR_r(x_1,\ldots,x_n),\]
whence $\nabla g_0(M\x)$ vanishes whenever $x_1,\ldots,x_r$ and
$R_1,\ldots,R_r$ all vanish.  This produces a subset of the singular
locus having projective dimension at least $n-1-2r\ge n-9$.  Hence if
$s(g_0)<n-9$, and in particular if $s(g_0)=-1$ and $n\ge 10$,
then there will be non-singular $p$-adic points whenever
the Congruence Condition holds.

We turn now to the singular integral $\mathfrak{I}(g;P)$.  
It follows from Lemma \ref{l8} that
\[\int_{-1}^{1}I(z;\ma{0})\d z=\int_{-P^{-11/4}}^{P^{-11/4}}I(z;\ma{0})\d z
+O(P^{7(n-3)/8}(\log P)^{7n}).\]
Moreover $g(\x)=g_0(\x)+O(|\x|^2)+O(1)$, whence 
\[e(zg(\x))=e(zg_0(\x))+O(|z|.|\x|^2)+O(|z|).  \]
It therefore follows that
\[I(z;\ma{0})=J(z)+O\Big(\int w(\x)|z|(|\x|^2+1)\d\x\Big)
=J(z)+O(|z|P^{n+2}),\]
where
\[J(z):=\int w(\x)e(z g_0(\x))\d\x.\] 
We now see that
\[\int_{-1}^{1}I(z;\ma{0})\d z=\int_{-P^{-11/4}}^{P^{-11/4}}J(z)\d z
+O(P^{7(n-3)/8}(\log P)^{7n})+O(P^{n-7/2}).\]
However Lemma \ref{l8} also applies to $J(z)$, whence
\[\int_{-P^{-11/4}}^{P^{-11/4}}J(z)\d z=\int_{-1}^{1}J(z)\d z
+O(P^{7(n-3)/8}(\log P)^{7n}).\]
Finally, from \cite[(10.3) and (10.4)]{hb-10} we see that
\[P^{n-3}(\log P)^{2-2n}\ll \int_{-1}^{1}J(z)dz\ll
P^{n-3}(\log P)^{2-2n}.\]
This therefore establishes \eqref{eq:n-3},
providing that $n\ge 4$, which thereby completes the proof of
Theorem \ref{main1'}.

\section{Proof of Theorem \ref{main1}: 
hyperplane sections}\label{s:slice}

It remains to prove Theorem \ref{main1}, which will be achieved by
induction on $s=s(g_0)$.  The base case will be $s=-1$, which follows
from Theorem \ref{main1'}.  For the induction we will find a
non-degenerate affine
hyperplane section of $g=0$ which again satisfies the Congruence
Condition, and for which $s$ is reduced by 1.

We begin by applying Bertini's Theorem (see Harris \cite[Theorem
17.16]{harris}, for example) to show that, for a generic vector
$\ma{a}$, the singular locus of the projective hyperplane section
$g_0(\x)=\ma{a}.\x=0$ has dimension $s-1$.  Similarly by Harris
\cite[Proposition 18.10]{harris}, for generic $\ma{a}$ the
intersection will be non-degenerate.  Thus there is a non-zero
form $f$ say, such that the dimension is $s-1$, and the intersection
is non-degenerate, whenever
$f(\ma{a})\not=0$.  Choose $\ma{a}$ to be any primitive integer vector such
that $f(\ma{a})\not=0$, whence the singular locus of 
$g_0(\x)=\ma{a}.\x=0$ will have dimension $s-1$.  The affine
hyperplane section we seek will then take the form $\ma{a}.\x=c$ for
a suitably chosen integer $c$.  We can find a matrix
$M\in {\rm SL}_n(\Z)$ whose first row is $\ma{a}$.  If we then write
$g$ in terms of $\y:=M\x$, by setting $h(\y)=g(M^{-1}\y)$, it follows
that the singular locus of $h_0(\y)=y_1=0$ will have dimension
$s-1$. Thus, irrespective of the value of $c$, if we put
$h^{(c)}(y_2,\ldots,y_n)=h(c,y_2,\ldots,y_n)$ then $s(h^{(c)}_0)=s-1$.
It is clear that distinct integer solutions $(y_2,\ldots,y_n)$ of
$h^{(c)}(y_2,\ldots,y_n)=0$, for some value of $c$, produce distinct 
solutions of $g(\x)=0$.

To complete our induction it therefore suffices to show that there is
an integer $c$ for which $h^{(c)}$ satisfies the Congruence Condition.
It will be convenient to use the notation $\u=(y_2,\ldots,y_n)$.
We begin by proving the following result.

\begin{lem}\label{bigp}
Suppose $h(y_1,\ldots,y_n)\in\Z[x_1,\ldots,x_n]$ is a cubic polynomial 
with $n\ge 4+s(h_0)$.  
Then there is an integer $p(h)$ depending on
$h$ such that for every prime $p\ge p(h)$ and every integer $c$ the congruence
\[h^{(c)}(\u)\equiv 0\mod{p}\]
has a non-singular solution modulo $p$.
\end{lem}

By Hensel's lemma, once we have a non-singular solution modulo $p$ we
will have solutions modulo $p^k$ for every $k$.

For the proof we define
$$
H^{(c)}(t,\u):=t^3h^{(c)}(t^{-1}\u),\quad
H_1(\u):=h^{(c)}_0(\u)=H^{(c)}(0,\u). 
$$
It will be important to observe that 
$H_1$ is independent of $c$, and that, by construction,
$s(H_1)=s(h_0)-1$. 
We proceed to estimate the number $N$, say, of solutions
to the congruence $h^{(c)}(\u)\equiv 0\tmod{p}$.  We have
$N=(N_1-N_2)/(p-1)$, where $N_1$ counts solutions of
\[H^{(c)}(t,\u)\equiv 0\mod{p} \]
and $N_2$ counts solutions of
$H_1(\u)\equiv 0\tmod{p}$.
For a form $F$ we shall write $s_p(F)$ to
denote the dimension of the singular locus of $F=0$ over $\ov{\F_p}$.  
If $p$ is sufficiently large then $s_p(H_1)=s(H_1)$, where 
``sufficiently large'' will be independent of $c$, since $H_1$ is
independent of $c$.  Taking hyperplane sections we can change the 
dimension of the singular locus by at most one, whence $s_p(H^{(c)})\le
1+s_p(H_1)=s(h_0)$ for large enough $p$.  Thus Lemma \ref{lemdel} yields
$N_1=p^{n-1}+O_n(p^{(n+1+s(h_0))/2})$ and 
$N_2=p^{n-2}+O_n(p^{(n-1+s(h_0))/2})$. 
Since $s(h_0)\le n-4$, by the hypothesis for Lemma \ref{bigp}, these 
bounds are enough to ensure that $N\gg_n p^{n-2}$ for large enough 
$p$.

To complete our treatment of ``large'' primes we estimate the number,
$S$ say, of singular solutions to $h^{(c)}(\u)\equiv 0\tmod{p}$.
Clearly $S\le S_1/(p-1)$ where $S_1$ is the number of solutions to
\[H^{(c)}=\frac{\partial H^{(c)}}{\partial u_1}=\ldots
=\frac{\partial H^{(c)}}{\partial u_{n-1}}=0\]
in $\F_p$.  Suppose these equations define a variety $V$ say in
projective space, and consider the variety $W$ defined by $\partial
H^{(c)}/\partial t=0$. 
Clearly $W$ has codimension at most 1, and
$V\cap W$ is the singular locus of $H^{(c)}$.  Thus $\dim(V\cap
W)=s_p(H^{(c)})\le s(h_0)$ for sufficiently large primes, 
as noted above, and hence $\dim(V)\le s(h_0)+1$.  It follows that
$S_1\ll_n p^{s(h_0)+2}$  
for large enough primes $p$, and hence that $S\ll_n p^{s(h_0)+1}$. 
Since $s(h_0)\leq n-4$ and $N\gg_n p^{n-2}$, 
we conclude that $N>S$ for large enough $p$, whence
$h^{(c)}(\u)\equiv 0\tmod{p}$ has a non-singular solution, as claimed.
This completes the proof of Lemma \ref{bigp}.

We now know that $h^{(c)}$ satisfies the Congruence Condition for 
$p\ge p(h)$ for every integer $c$.  To complete our argument we
proceed to choose $c$ so that the condition is satisfied for the
remaining small primes.  Now we saw in \S \ref{endgame} that if the
Congruence Condition holds for $g$ then there will in fact be
non-singular $p$-adic solutions, if $s=s(g_0)<n-9$.  
Thus, under the
hypotheses of Theorem \ref{main1}, we may assume that there is a
non-singular $p$-adic solution of $g=0$ 
for each $p$.  It then follows that $h(\y)=0$
has a non-singular $p$-adic solution $\y_0$, say.  We will need to know
that there must be a solution with $\nabla' h(\y)\not=\ma{0}$, where
\[\nabla' h(\y):=\Big(\frac{\partial h}{\partial y_2},\ldots,
\frac{\partial h}{\partial y_n}\Big).\]
This is the content of the following result.

\begin{lem}\label{last}
Suppose $h(y_1,\ldots,y_n)\in\Z_p[x_1,\ldots,x_n]$ is a cubic polynomial 
with $h_0$ absolutely irreducible.  
If $h(\y)=0$ has a non-singular solution in
$\Z_p$ then there is a solution with $\nabla' h(\y)\not=\ma{0}$.
\end{lem}

We argue by contradiction.  Suppose that $\nabla' h(\y_0)=\ma{0}$ for
every non-singular $p$-adic solution $h(\y_0)=0$.  
Let $\w$ be a $p$-adic integer
vector with $w_1=0$.  Then if $|\w|_p$ is sufficiently small,
Hensel's Lemma shows that there is a non-singular solution $h(\y_0+\w+\z)=0$
with $|\z|_p\ll |\w|_p^2$.  We supposedly have $\nabla' 
h(\y_0+\w+\z)=\ma{0}$.  However if $M=M(\y_0)$ 
is the matrix of second derivatives of $h$ at $\y_0$, then
\[\nabla h(\y_0+\w+\z)=\nabla h(\y_0)+M\w+O(|\w|_p^2).\]
(By this, we mean that one can replace the error term by a vector
whose $p$-adic norm is $O(|\w|_p^2)$.)
Since $\nabla' h(\y_0+\w+\z)=\nabla' h(\y_0)=\ma{0}$ we deduce that
$(M\w)_i=O(|\w|_p^2)$ for $2\le i\le n$, whenever $w_1=0$ and 
$|\w|_p$ is small enough.  Since $\w$ is arbitrary subject to these
restrictions it follows that $M_{ij}(\y_0)=0$ whenever
$2\le i,j\le n$.

However $\y_0$ was an arbitrary non-singular solution of $h(\y_0)=0$.
It therefore follows that $M_{ij}(\y)=0$ for $2\le i,j\le n$, for every 
non-singular solution of $h(\y)=0$.  We may therefore repeat our
argument.  Since $h(\y_0+\w+\z)=0$, we deduce that 
$M_{ij}(\y_0+\w+\z)=0$.  This time we have
\[M_{ij}(\y_0+\w+\z)=M_{ij}(\y_0)+6\sum_{k=1}^n c_{ijk}w_k+O(|\w|_p^2),\]
if
\[h_0(\y)=\sum_{i,j,k=1}^nc_{ijk}y_iy_jy_k\]
with symmetric coefficients $c_{ijk}$.  Arguing as before we deduce that
\[\sum_{k=1}^n c_{ijk}w_k=O(|\w|_p^2), \quad (2\le i,j\le n), \]
whenever $w_1=0$ and $|\w|_p$ is small enough. This
allows us to conclude that $c_{ijk}=0$ for $2\le i,j,k\le n$.
It therefore follows that
$y_1$ divides $h_0(\y)$ identically.  We have finally reached a
contradiction, and the lemma follows.

In our situation, if $h_0$ were reducible, we would have $s(h_0)\ge
n-3$, which is contrary to hypothesis.  
Lemma \ref{last} therefore implies that for 
each $p< p(h)$ we can find a
$p$-adic integer vector $\y^{(p)}$ with $h(\y^{(p)})=0$ and $\nabla'
h(\y^{(p)})\not=\ma{0}$.  Suppose that the exponent $k(p)$ satisfies 
$p^{k(p)}\mid \nabla' h(\y^{(p)})$ but $p^{k(p)+1}\nmid \nabla' h(\y^{(p)})$, 
and choose a vector $\z^{(p)}\in\Z^n$ with
$\z^{(p)}\equiv\y^{(p)}\tmod{p^{2k(p)+1}}$.  We define 
\[\u^{(p)}= \big(z^{(p)}_2,\ldots,z^{(p)}_n\big).\]
Finally let $c$ satisfy $c\equiv
z^{(p)}_1\tmod{p^{2k(p)+1}}$ for every prime $p< p(h)$.  Such a $c$
exists, by the Chinese Remainder Theorem.  Then
\begin{align*}
h^{(c)}(\u^{(p)})=h(c,\u^{(p)})&\equiv h(z^{(p)}_1,\u^{(p)}) \mod{p^{2k(p)+1}}\\
&=h(\z^{(p)})\\
&\equiv h(\y^{(p)}) \mod{p^{2k(p)+1}}\\
&\equiv 0\mod{p^{2k(p)+1}},
\end{align*}
while
\[\nabla h^{(c)}(\u^{(p)})=\nabla h(c,\u^{(p)})\equiv \nabla' h(\z^{(p)})\equiv
\nabla' h(\y^{(p)})\not\equiv \ma{0} 
\mod{p^{k(p)+1}}.\]
It follows that the vector $\u^{(p)}$ can be lifted to a non-singular $p$-adic
solution of $h^{(c)}(\u)=0$.  This establishes the Congruence
Condition for $h^{(c)}$ for every prime $p<p(h)$, thereby completing
the proof of Theorem \ref{main1}.

\end{document}